\documentclass[twocolumn]{autart}    

\usepackage{graphicx}           
\usepackage{dcolumn}            
\usepackage{bm}                 

\usepackage{graphics}           
\usepackage{epsfig}             
\usepackage{times}              
\usepackage{amsmath}
\usepackage{amssymb}
\usepackage{extarrows}
\usepackage{amsfonts}
\usepackage{fancyhdr}
\usepackage{color}
\usepackage{subfigure}
\usepackage{enumerate}

\allowdisplaybreaks[4] 

\newtheorem{theorem}{Theorem}
\newtheorem{lemma}[theorem]{Lemma}

\newtheorem{proposition}[theorem]{Proposition}

\newtheorem{definition}{Definition}

\newtheorem{remark}{Remark}

\newtheorem{assumption}{Assumption}

\begin{document}

\begin{frontmatter}

\title{\textcolor{black}{Input-to-State} Stability with Respect to Boundary Disturbances for a Class of Semi-linear Parabolic Equations} 

\thanks[footnoteinfo]{Corresponding author: G. Zhu}

\author{Jun Zheng$^{1,2}$}\ead{zhengjun2014@aliyun.com},
\author{Guchuan Zhu$^{3}$}\ead{guchuan.zhu@polymtl.ca}

\address{$^{1}$Department of Basic Course, Southwest Jiaotong University\\
        Emeishan, Sichuan, P. R. of China 614202\\
        $^{2}${School of Mathematics, Southwest Jiaotong University\\
        Chengdu, Sichuan, P. R. of China 611756}\\
        $^{3}$Department of Electrical Engineering, Polytechnique
        Montr\'{e}al, \\ P.O. Box 6079, Station Centre-Ville,
        Montreal, QC, Canada H3T 1J4}  

\begin{keyword}                            
Semi-linear parabolic equations, %
ISS, boundary disturbances, Lyapunov method. 
\end{keyword}                              

\begin{abstract}                           
This paper studies the \textcolor{black}{input-to-state stability} (ISS) properties \textcolor{black}{based on the method of Lyapunov functionals} for a class of semi-linear parabolic partial differential equations (PDEs) \textcolor{black}{with respect to boundary disturbances}. In order to avoid the appearance of time derivatives of the disturbances in ISS estimates, some technical inequalities are first developed, which \textcolor{black}{allow} directly dealing with the boundary conditions \textcolor{black}{and establishing the ISS based on the method of Lyapunov functionals}. The well-posedness analysis of the considered problem is carried out and the \textcolor{black}{conditions for ISS are derived}. Two examples are used to illustrate the application of the developed result.
\end{abstract}
\end{frontmatter}
\section{Introduction}\label{Sec: Introduction}
In the past few years, there has been a considerable effort devoted to extending the input-to-state stability (ISS) theory, which was originally introduced by Sontag for finite-dimensional nonlinear systems \cite{Sontag:1989,Sontag:1990}, to infinite dimensional systems governed by partial differential equations (PDEs). In particular, significant progresses on the establishment of ISS properties with respect to disturbances for different PDEs have been reported in the recent literature   \cite{Argomedo:2013,Argomedo:2012,Dashkovskiy:2010,Dashkovskiy:2013,Dashkovskiy:2013b,jacob2016input,Karafyllis:2014,Karafyllis:2016,Karafyllis:2016a,karafyllis2017siam,Karafyllis:2017,Logemann:2013,Mazenc:2011,Mironchenko:2014b,Mironchenko:2014a,Prieur:2012,Tanwani:2017}.

It is noticed that the majority of the existing work dealt with disturbances distributed over the domain for which the method of Lyapunov functionals is shown to be a well-suited tool. However, difficulties may be encountered when considering disturbances acting on the boundaries. This is mainly due to the fact that the latter case usually leads to a formulation involving unbounded operators, which may be an obstacle for the construction of Lyapunov functionals as explained in \cite{Karafyllis:2016,Karafyllis:2016a,karafyllis2017siam}.
It is shown in \cite{Jacob:2016,jacob2016input} that for a class of linear PDEs, the exponential stability plus a certain admissibility implies the ISS with respect to boundary disturbances. \textcolor{black}{However, it may be difficult to characterize the admissibility for nonlinear PDEs.} To avoid dealing with unbounded operators, it is proposed in \cite{Argomedo:2012} to transform the boundary disturbance to a distributed one, which allows for the application of the well-established tools, in particular the method of Lyapunov functionals. However, as pointed out in \cite{Karafyllis:2016,Karafyllis:2016a,karafyllis2017siam} \textcolor{black}{the result given in \cite{Argomedo:2012} may} end up with ISS estimates expressed by boundary disturbances and their time derivatives, which is not strictly in the original form of ISS formulation. To resolve this concern, it is proposed in \cite{Karafyllis:2016,Karafyllis:2016a,karafyllis2017siam} to derive the ISS property directly from the estimates of the solution to the considered PDEs by using eigenfunction expansions or finite-difference schemes. \textcolor{black}{An advantage of these methods is that they can be applied \textcolor{black}{to a wide range of linear and nonlinear PDEs. Whereas,} these methods may involve heavy computations. In a recent work \cite{Mironchenko:2017}, a new method based monotonicity has been introduced for studying the ISS of nonlinear parabolic equations with boundary disturbances. As an application of this method, the ISS properties in $L^p$-norm ($p>2$) for some linear parabolic equations with Dirichlet boundary disturbances have been established.}
Nevertheless, it is still of great interest to investigate the applicability of the well-established method of Lyapunov functionals to the establishment of ISS properties with respect to boundary disturbances \textcolor{black}{for nonlinear PDEs, including those investigated recently in \cite{Karafyllis:2016,Karafyllis:2016a,karafyllis2017siam,Mironchenko:2017}}. This motivates the present work.

The aim of this work is to establish the aforementioned ISS property for a class of semi-linear parabolic PDEs \textcolor{black}{with Robin (or Neumann) boundary conditions based on the method of Lyapunov functionals}. To achieve this objective, we have developed first in Section~\ref{Sec: preliminary} some technical \textcolor{black}{inequalities }(Lemma~\ref{Lemma 3} and Lemma~\ref{Lemma 2}) {that establish some relationships between the value of a real-valued $C^1$-function at any point and its norms. This is a key feature that allows dealing directly with the boundary conditions and avoiding the appearance of time derivatives of the disturbance in ISS estimates.} The well-posedness of the problem described in Section~\ref{Sec: Problem formulation} is addressed in Section~\ref{Sec: Well-posedness}. A quite standard Lyapunov functional \cite{Mazenc:2011} is then used in Section~\ref{Sec: ISS} to establish the ISS estimates of the solutions with respect to in-domain and boundary disturbances. Finally, the ISS analysis of two parabolic PDEs are given in Section~\ref{Sec: Examples} to illustrate the proposed method. The main contribution of the present work is the derivation of the ISS property of the considered PDEs from a Lyapunov functional using the developed techniques that can be useful in the study of other types of PDEs.

\textbf{Notation.} In this paper, $\mathbb{R}_+$ denotes the set of positive real numbers and $\mathbb{R}_{\geq 0} :=  0\cup\mathbb{R}_+$.
$L^2(a,b)=\{u:(a,b)\rightarrow\mathbb{R}
|\ \int_{a}^bu^2(x)\text{d}x<+\infty\}$, which is a Hilbert space endowed with the natural product $\langle u,v\rangle= \int_{a}^bu(x)v(x)\text{d}x$.
$H^2(a,b)=\{u:(a,b)\rightarrow\mathbb{R}
|\ u\in L^2(a,b)$ and the derivatives of first order and second order $u_x,u_{xx}\in L^2(a,b)\}$.
$C(\mathbb{R}_{\geq 0};\mathbb{R})=\{u:\mathbb{R}_{\geq 0}\rightarrow\mathbb{R}|$ $u$ is continuous on $\mathbb{R}_{\geq 0}\}$.
$C^1([a,b];\mathbb{R})=\{u:[a,b]\rightarrow\mathbb{R}|\ u$ and $u_x$ are continuous on $[a,b]\}$.
$C^2(\mathbb{R}_{\geq 0};\mathbb{R})=\{u:\mathbb{R}_{\geq 0}\rightarrow\mathbb{R}|\ u$, $u_x$ and $u_{xx}$ are continuous on $\mathbb{R}_{\geq 0}\}$.
$C^0(\mathbb{R}_{\geq 0}\times(a,b);\mathbb{H})=\{u:\mathbb{R}_{\geq 0}\times(a,b)\rightarrow \mathbb{R}| \ u(t,\cdot)\rightarrow\mathbb{H},\ u$ is continuous (in $t$) on $\mathbb{R}_{\geq 0}\}$, where $\mathbb{H} $ is some \textcolor{black}{function} space.
$C^1(\mathbb{R}_{+}\times(a,b);\mathbb{H})=\{u:\mathbb{R}_{+}\times(a,b)\rightarrow \mathbb{R}| \ u(t,\cdot)\rightarrow\mathbb{H},u_t(t,\cdot)\rightarrow\mathbb{H},\ u$ and $u_t$ are continuous (in $t$) on $\mathbb{R}_{+}\}$, where $\mathbb{H} $ is some \textcolor{black}{function} space. \textcolor{black}{Let $\mathcal {K}=\{\gamma : \mathbb{R}_{\geq 0} \rightarrow \mathbb{R}_{\geq 0}|\ \gamma(0)=0,\gamma$ is continuous, strictly increasing$\}$; $ \mathcal {K}_{\infty}=\{\theta \in \mathcal {K}|\ \lim\limits_{s\rightarrow\infty}\theta(s)=\infty\}$; $ \mathcal {L}=\{\gamma : \mathbb{R}_{\geq 0}\rightarrow \mathbb{R}_{\geq 0}|\ \gamma$ is continuous, strictly decreasing, $\lim\limits_{s\rightarrow\infty}\gamma(s)=0\}$; $ \mathcal {K}\mathcal {L}=\{\beta : \mathbb{R}_{\geq 0}\times \mathbb{R}_{\geq 0}\rightarrow \mathbb{R}_{\geq 0}|\ \beta(\cdot,t)\in \mathcal {K}, \forall t \in \mathbb{R}_{\geq 0}$, and $\beta(s,\cdot)\in \mathcal {L}, \forall s \in \mathbb{R}_{\geq 0}\}$.}

Throughout this paper, we always denote $\|u\|_{L^2(a,b)} $, or $\|u\|_{L^2(0,1)}$, by $\|u\|$ for notational simplicity.
\section{Problem setting and preliminaries}\label{Sec: preliminary}
\subsection{Problem setting}\label{Sec: Problem formulation}
We consider the following semi-linear 1-$D$ parabolic equation
\begin{align}\label{Equ. 1}
 u_t-\mu u_{xx}=f(t,x,u,u_x) \ \ \ \ \text{in}\ \  \mathbb{R}_{\geq 0}\times (0,1)
\end{align}
with the boundary and initial conditions
\begin{subequations}\label{B}
\begin{align}
 &a_1u(t,1)+a_2u_x(t,1)=0,\label{B. 1}\\
 &b_1u(t,0)+b_2u_x(t,0)=d_1(t),\label{B. 2}\\
  &u(0,x)=u_0(x),\label{B. 3}
\end{align}
\end{subequations}
where $d_1(t)$ is the disturbance acting on the boundary, $a_1,a_2,b_1$, $b_2$ {are nonnegative constants and $\mu$ is a positive constant.} In \eqref{Equ. 1}, \textcolor{black}{for the function $f:\mathbb{R}_{\geq 0}\times (0,1)\times\mathbb{R}\times\mathbb{R}\rightarrow \mathbb{R}$, there exist a continuous function
$\rho:\mathbb{R}_{\geq 0}\times\mathbb{R}\rightarrow\mathbb{R}_{\geq 0}$, which is monotonously increasing in the second argument, a constant $\gamma \in [1,3)$, and a constant $\vartheta \in (0,1]$, such that for any $T\in \mathbb{R}_{+}$, there hold
\begin{subequations}\label{AF}
\begin{align}
&|f(t,x,u,p)|\leq \rho(t,|u|)(1+|p|^{\gamma}),\label{AF-1}\\
&|f(s,x,u,p)-f(t,x,u,q)|\notag\\
&\;\;\;\;\;\;\; \leq \rho(0,|u|)(1+|p|^{\gamma})|s-t|^{\vartheta},\label{AF-2'}\\
&|f(t,x,u,p)-f(t,x,u,q)|\notag\\
&\;\;\;\;\;\;\; \leq \rho(t,|u|)(1+|p|^{\gamma-1}+|q|^{\gamma-1})|p-q|,\label{AF-2}\\
&|f(t,x,u,p)-f(t,x,v,p)| \notag\\
&\;\;\;\;\;\;\; \leq \rho(t,|u|+|v|)(1+|p|^{\gamma})|u-v|,\label{AF-3}
\end{align}
\end{subequations}
for} a.e. $x\in (0,1)$ and all $ s,t\in [0,T) , u\in\mathbb{R},v\in\mathbb{R},p\in\mathbb{R}$.
\subsection{Preliminaries}\label{Sec: preliminary}
\vspace*{-6pt}
In the subsequent development, we employ extensively the following inequalities.
\vspace*{6pt}
\newline
\textbf{Young's inequality:} For real numbers $a\geq 0$, $b\geq 0$, and $\varepsilon>0$, there holds $ab\leq \frac{a^2}{2\varepsilon}+\frac{\varepsilon b^2}{2}$.
%
\vspace*{6pt}
\newline
\textbf{Gronwall's inequality:} Suppose that $y:\mathbb{R}_{\geq 0}\rightarrow\mathbb{R}_{\geq 0}$ is absolutely continuous on $[0,T]$ for any $T>0$ and satisfies for a.e. $t \geq 0$ the differential inequality
$
\frac{\text{d}y}{\text{d}t}(t)\leq g(t)y(t)+h(t),
$
where $g,h\in L^1([0,T];\mathbb{R})$ for any $T > 0$. Then for all $t \in \mathbb{R}_{\geq 0}$, there holds
\begin{equation*}
  y(t) \leq y(0)e^{\int_{0}^{t}g(s)\text{d}s}+\displaystyle\int_{0}^{t}h(s) e^{\int_{s}^{t}g(\tau)\text{d}\tau}\text{d}s.
\end{equation*}
The following inequalities will be used to deal with the items associated with boundary points. They are essential for establishing the ISS property with respect to boundary disturbances without invoking their time derivatives in \emph{a priori} estimates of the solution.
%
%
\begin{lemma}\label{Lemma 3}
\textcolor{black}{Suppose that $u\in C^{1}([a,b];\mathbb{R})$, then
\begin{align*}
u^2(c)\leq \frac{2}{b-a}\|u\|^2+(b-a)\|u_x\|^2,\ \ \ \ \forall c\in[a,b].
\end{align*}}
\end{lemma}
\vspace*{-24pt}
\begin{pf*}{Proof.}
\textcolor{black}{For each $c\in [a,b]$, let} $g(x)=\textcolor{black}{\int_{c}^xu^2_z(z)\text{d}z}$.
\textcolor{black}{ Note that $g_x(x)=u^2_x(x)$. By Cauchy-Schwarz
inequality, we have
\begin{align*}
\bigg(\int_{c}^xu_z(z)\text{d}z\bigg)^2&\leq \bigg|(x-c)\int_{c}^xu^2_z(z)\text{d}z\bigg|\notag\\
&=(x-c)\int_{c}^xu^2_z(z)\text{d}z.
\end{align*}
It follows
\begin{align}
&\int_{a}^b\bigg(\int_{c}^xu^2_z(z)\text{d}z\bigg)^2\text{d}x \leq \int_{a}^b(x-c) g(x)\text{d}x\notag\\
=&\bigg[\frac{(x-c)^2}{2}g(x)\bigg]\bigg|_{x=a}^{x=b}-\int_{a}^b\frac{(x-c)^2}{2} u^2_x(x) \text{d}x\notag\\
\leq &\frac{(b-c)^2}{2}\int_{c}^bu^2_x(x)\text{d}x- \frac{(a-c)^2}{2}\int_{c}^au^2_x(x)\text{d}x\notag\\
=&\frac{(b-c)^2}{2}\int_{c}^bu^2_x(x)\text{d}x+ \frac{(a-c)^2}{2}\int_{a}^cu^2_x(x)\text{d}x\notag\\
\leq &\frac{(b-a)^2}{2}\int_{a}^bu^2_x(x)\text{d}x.\label{6262}
\end{align}
Note that
\begin{align*}
u^2(c)\!=\!\bigg(\!\!u(x)+\!\!\textcolor{black}{\int_{x}^c\!\!u^2_z(z)\text{d}z}\bigg)^2\!\!\leq \! 2u^2(x)+2\bigg(\textcolor{black}{\int_{x}^c\!\!u^2_z(z)\text{d}z}\bigg)^2.
\end{align*}
Integrating over $[a,b]$ and noting \eqref{6262}, we get
\begin{align*}
u^2(c)(b-a)\leq 2\int_{a}^bu^2(x)\text{d}x+(b-a)^2\int_{a}^bu^2_x(x)\text{d}x,
\end{align*}
which yields the claimed result.} \hfill $\blacksquare$
\end{pf*}
\begin{lemma}\label{Lemma 2}
Suppose that $u\in C^{1}([a,b];\mathbb{R})$.
\begin{itemize}
\item[(i)] \textcolor{black}{ If $u(c_0)=0$ for some $c_0\in[a,b]$, there holds}
$\textcolor{black}{\|u\|^2\leq  \frac{(b-a)^2}{2}\|u_x\|^2.}$
\item[(ii)] For any {$c\in[a,b]$, there holds}
 $\|u\|^2\leq  2u^2(c)(b-a)+(b-a)^2\|u_x\|^2.$
%
\end{itemize}
\end{lemma}
\begin{pf*}{Proof.}
\textcolor{black}{Note that for any $w\in C^1([0,1]) $, there holds \cite{Krstic:2008b}}
\begin{align*}
\textcolor{black}{\|w\|^2_{L^2(0,1)}\leq w^2(i)+\frac{1}{2}\|w_x\|^2_{L^2(0,1)},\ i=0,1.}
\end{align*}
\textcolor{black}{For any $c\in [a,b]$, let $v(x)=u(c-(c-a)x)$. Then we get}
\textcolor{black}{\begin{align*}
\|u\|^2_{L^2(a,c)}&=(c-a)\|v\|^2_{L^2(0,1)}\\
&\leq  (c-a) v^2(0)+ \frac{c-a}{2}\|v_x\|^2_{L^2(0,1)}\\
&=(c-a) u^2(c)+\frac{(c-a)^2}{2}\|u_x\|^2_{L^2(a,c)}.
\end{align*}}
\textcolor{black}{Similarly, one may get
\begin{align*}
\|u\|^2_{L^2(c,b)}\leq (b-c) u^2(c)+\frac{(b-c)^2}{2}\|u_x\|^2_{L^2(c,b)}.
\end{align*}
Finally, one has
\begin{align*}
\|u\|^2_{L^2(a,b)}&=  \|u\|^2_{L^2(a,c)}+\|u\|^2_{L^2(c,b)}\\
&\leq (b-a) u^2(c)+\frac{(b-a)^2}{2}\|u\|^2_{L^2(a,b)}.
\end{align*}}
\textcolor{black}{Then, we can conclude that (i) and (ii) hold true.}
\hfill $\blacksquare$
\end{pf*}
\section{Well-posedness Analysis}\label{Sec: Well-posedness}
Consider first the solution to \eqref{Equ. 1} with disturbance free boundary conditions:
\begin{subequations}\label{Equ. 2}
\begin{align}
 &a_1u(t,1)+a_2u_x(t,1)=0,\\
 & b_1u(t,0)+b_2u_x(t,0)=0,\\
 &u(0,x)=\textcolor{black}{u_0(x)}.
\end{align}
\end{subequations}
In this section, we always assume that
\begin{align*}
&{(a_1+a_2)b_1\neq a_1b_2},\ d_{1}\in C^{2}(\mathbb{R}_{\geq 0};\mathbb{R}),\label{AB} \\
&u_0\in \mathbb{H}^2_{(0)}:=\{u\in H^2(0,1); a_1u(1)+a_2u_x(1)= 0, \notag\\
&\;\;\;\;\;\;\;\;\;\;\;\;\;\;\;\;\;\;\;\;\;\;\;\;\;\;\;\;\;\; b_1u(0)+b_2u_x(0)=0\}.\notag
\end{align*}
Moreover, we make the following assumptions.
When $a_2b_2= 0$, we always assume that
\begin{align*}
\frac{a_1}{a_2}&\geq -\frac{1}{2}, &\text{if}\ \ a_2\neq 0\ \ \text{and}\ \ b_2=0, \\
\frac{b_1}{b_2}&\leq \frac{1}{2},  &\text{if}\ \ a_2=0\ \ \text{and}\ \ b_2\neq0.
\end{align*}
%
%
%
When $a_2b_2\neq 0$, we always assume that there \textcolor{black}{exist} $A_1, A_2\in \mathbb{R}_{\geq 0}$ satisfying $A_1+A_2=1 $ such that
\begin{equation}\label{+ab1}
 \frac{a_1}{a_2} \geq 2A_2,\ \frac{b_1}{b_2} \leq A_1,\ A_2-2A_1 \geq 0,
\end{equation}
or, there \textcolor{black}{exist} $B_1,B_2\in \mathbb{R}_{\geq 0}$ satisfying $B_1+B_2=1 $ such that
\begin{align}
  \frac{a_1}{a_2} &\geq -B_2,\ \frac{b_1}{b_2}\leq -2B_1,\ B_1-2B_2\geq 0.\label{+ab2}
\end{align}
%
\begin{remark}\label{Remark 1}
Under the above assumptions, we always have
$a_1^2+a_2^2>0$, $b_1^2+b_2^2>0$ and $a_1^2+b_1^2>0$.
\end{remark}
\begin{proposition}\label{Proposition 1} Assume that $u_1\in \mathbb{H}^2_{(0)}$. Then there exists a unique solution $u\in C^{0}(\mathbb{R}_{\geq 0}\times (0,1);\mathbb{H}^2_{(0)})\cap C^{1}(\mathbb{R}_{+}\times (0,1);L^2(0,1))$ to \eqref{Equ. 1} with boundary-initial conditions \eqref{Equ. 2}.
\end{proposition}
\begin{pf*}{Proof.} \textcolor{black}{Let $D(\mathcal{A})=\mathbb{H}^2_{(0)}$ with the operator $\mathcal{A}: D(\mathcal{A})\rightarrow L^2(0,1)$ defined as $\mathcal{A}u=\mu u_{xx}$. For $\alpha\in [0,1]$, let $\mathcal{H}_{\alpha}=\mu (-\mathcal{A})^{\alpha} $ and $\|u\|_{\alpha}=\|(-\mathcal{A})^{\alpha}u\|$. Let $F(t,u)(x)=f(t,x,u,u_x)$. Then \eqref{Equ. 2} can be expressed by an abstract evolutionary equation $\frac{\text{d}u}{\text{d}t}=\mathcal{A}u+F(t,u)$. The proof is based on the theory of Lipschitz perturbations of linear evolution equations \cite[Theorem 12, \S4.3]{Kuttler:2011} (see also \cite[\S6.3, Chap.~6]{Pazy:1983}),} \textcolor{black}{which consists in two steps: first to prove that $\mathcal{A}$ is the infinitesimal generator of a $C_0$-semigroup of contractions on $L^2(0,1)$; and second to prove that $F(t,u) $ satisfies local H\"{o}lder condition, i.e., for $F:\mathbb{R}_{\geq 0}\times U\rightarrow L^2(0,1)$, where $U$ is an open subset of $\mathcal{H}_{\alpha}$, for every $(t,u)\in U$, there is a neighborhood $V\subset U$ and constants $L\geq 0$, $0<\vartheta\leq 1$ (see, e.g., \cite[Assumption (NONLIN), \S4.3]{Kuttler:2011}, \cite[Assumption (F), \S6.3, Chap.~6]{Pazy:1983}) such that
\begin{align}
\|F(t_1,u_1)-F(t_2,u_2)\|
\leq &L(|t_1 -t_2|^{\vartheta}+\|u_1-u_2\|_{\alpha}),\notag\\
& \forall (t_i,u_i)\in V,i=1,2.\label{+9251}
\end{align}
First, since $\mathcal{A}$ is a densely defined closed linear operator and self-adjoint, it suffices to prove that $\mathcal{A}$ is dissipative. Then, the claim that $\mathcal{A}$ generates a $C_0$-semigroup follows} from Lumer-Phillips theorem (see \cite[Corollary 4.4, \S1.4, Chap.~1]{Pazy:1983}, \cite[Theorem 6.1.8]{Jacob:2012}).
Indeed, due to
\begin{align*}
\frac{1}{\mu}\langle \mathcal{A}u,u\rangle =& \int_{0}^1 u_{xx}u\text{d}x \\
                                 =&  u_{x}(t,1)u(t,1)- u_{x}(t,0)u(t,0)- \|u_x\|^2,
\end{align*}
we may argue for four cases.
\vspace*{6pt}
\newline
(i) $b_2=a_2=0$. In this case, $u(t,1)=u(t,0)=0$. $\mathcal{A}$ is obviously dissipative.
\vspace*{6pt}
\newline
(ii) $b_2=0,a_2\neq 0$. In this case, $u(t,0)=0$. It follows
\begin{align*}
\frac{1}{\mu}\langle \mathcal{A}u,u\rangle =\int_{0}^1
u_{xx}u\text{d}x= -\frac{a_1}{a_2}u^2(t,1)- \|u_x\|^2.
\end{align*}
For ${\frac{a_1}{a_2}\geq -\frac{1}{2}}$, $ -\frac{a_1}{a_2}u^2(t,1)\leq \frac{1}{2}u^2(t,1)$. By Lemma \ref{Lemma 3} and Lemma \ref{Lemma 2}, we get
$
u^2(t,1)\leq 2\|u\|^2+\|u_x\|^2\leq 2\|u_x\|^2.
$
Then
\begin{align*}
\frac{1}{\mu}\langle \mathcal{A}u,u\rangle =&-\frac{a_1}{a_2}u^2(t,1)- \|u_x\|^2 \\
 \leq & \frac{1}{2}\times 2\|u_x\|^2- \|u_x\|^2= 0.
\end{align*}
(iii) $b_2\neq0,a_2= 0$. In this case, $u(t,1)=0$.
Arguing as in (ii), for {$\frac{b_1}{b_2}\leq \frac{1}{2}$}, one may get
$
\frac{1}{\mu}\langle \mathcal{A}u,u\rangle
\leq 0.
$
\vspace*{6pt}
\newline
 (iv) $b_2\neq0,a_2\neq  0$.
For $\frac{a_1}{a_2}\geq 2A_2,\ \frac{b_1}{b_2}\leq A_1,\ A_2-2A_1\geq 0$, note that by Lemma \ref{Lemma 3} and Lemma \ref{Lemma 2} there hold
\begin{align}\label{62812}
\begin{split}
\|u_x\|^2\geq u^2(t,0)-2\|u\|^2,\\
\|u_x\|^2\geq \|u\|^2-2u^2(t,1).
\end{split}
\end{align}
Then we get
\begin{align*}
&\frac{1}{\mu}\langle \mathcal{A}u,u\rangle \notag\\
=& -\frac{a_1}{a_2}u^2(t,1)+\frac{b_1}{b_2}u^2(t,0)-A_1\|u_x\|^2-A_2\|u_x\|^2\notag\\
\leq & -\frac{a_1}{a_2}u^2(t,1)+\frac{b_1}{b_2}u^2(t,0)-A_2( \|u\|^2-2u^2(t,1)) \notag\\
& \;\;\; -A_1( u^2(t,0)-2\|u\|^2)\notag\\
=&\bigg( 2A_2-\frac{a_1}{a_2}\bigg)u^2(t,1)+\bigg( \frac{b_1}{b_2}-A_1\bigg)u^2(t,0)\notag\\
& \;\;\; +(2A_1-A_2)\|u\|^2 \leq 0.
\end{align*}
Similarly, for $\frac{a_1}{a_2}\geq -B_2,\ \frac{b_1}{b_2}\leq -2B_1,\ B_1-2B_2\geq 0$, 
\textcolor{black}{one may get
$\frac{1}{\mu}\langle \mathcal{A}u,u\rangle\leq 0.$}
Thus, $\mathcal{A}$ is a dissipative operator.

\textcolor{black}{The second step of proof can be proceeded in the same way as in \cite[Proposition 7, \S4.4]{Kuttler:2011}. First, since $\mathcal{A} $ is a Sturm-Liouville operator \cite{Boyce:1997,Naylor:1982}, all eigenvalues of $\mathcal{A} $ are real, and form an infinite, increasing sequence $0>\lambda_1>\lambda_2>\cdots>\lambda_n>\cdots$ with $\lim\limits_{n\rightarrow\infty}\lambda_n=-\infty $. Corresponding to each $\lambda_n\in \mathbb{R}, n=1,2,\ldots$, there is exactly one eigenfunction $ \varphi_n \in D(\mathcal{A})\cap C^2([0,1]) $ satisfying $ \mathcal{A}\varphi_n=\lambda_n\varphi_n$. The eigenfunctions form an orthonormal basis of $ L^2(0,1)$. Second, one may proceed exactly as in  \cite[\S4.3,\S4.4]{Kuttler:2011} to show that the norm $ \|u\|+\|u\|_{\alpha}$ on $ \mathcal{H}_{\alpha}$ is equivalent to the norm $\|u\|_{\alpha}$ and $ \mathcal{H}_{\alpha}\subset W^{1,2\gamma}(0,1)\cap L^{\infty}(0,1)$ for $ \max\{\frac{3}{4},\frac{5\gamma-3}{4\gamma}\}<\alpha<1$. Furthermore, Theorem 12 in \cite[\S4.3]{Kuttler:2011} holds. Then proceeding exactly as in (4.17)-(4.20) in \cite[Theorem 4.4, \S8.4, Chap.~8]{Pazy:1983}, one may verify that $F(t,u)$ satisfies \eqref{+9251}.}
\vspace*{6pt}
\newline
\textcolor{black}{ Finally, Theorem 12 in \cite[\S4.3]{Kuttler:2011} guarantees the existence of a unique classical solution.}\hfill $\blacksquare$

\end{pf*}
\begin{remark}
Conditions on the constants $a_i, b_i (i=1,2)$ are only required in order to guarantee exponential stability
of a semigroup in the proof of Proposition \ref{Proposition 1}.
\end{remark}
\begin{theorem}\label{Theorem 1}
There exists a unique solution $u\in C^{0}(\mathbb{R}_{\geq 0}\times (0,1);\mathbb{H}^2_{(0)})\cap C^{1}(\mathbb{R}_{+}\times (0,1);L^2(0,1))$ of \eqref{Equ. 1} satisfying \eqref{B. 1}, \eqref{B. 2}, and \eqref{B. 3}.
\end{theorem}
\begin{pf*}{Proof.}
Consider first the case where $b_1^2+b_2^2=1$. Let $g(x)=b_1+b_2x+c_1x^2+c_2x^3$, where $c_1, c_2\in \mathbb{R}$ satisfy $(a_1+2a_2)c_1+(a_1+3a_2)c_2=-a_1b_1-(a_1+a_2)b_2$. The existence of
$c_1, c_2$ is guaranteed by $a_1^2+a_2^2\neq 0$. One may check that
$a_1g(1)+a_2g_x(1)=b_1g(0)+b_2g_x(0)=0$ due to $b_1^2+b_2^2=1$.
\newline
Let $\tilde{f}(t,x,v,p)=d_{1t}(t)g(x)+\mu d_1(t)g_{xx}(x)+f(t,x,v+d_1(t)g(x),p+d_1(t)g_x(x))$. Consider the following equation
\begin{subequations}\label{Equ. 4}
\begin{align}
&v_t-\mu v_{xx}=\tilde{f}(t,x,v,v_x),\\
 &a_1v(t,1)+a_2v_x(t,1)=0,\\
 & b_1v(t,0)+b_2v_x(t,0)=0,\\
 &v(0,x)=v_0(x),
\end{align}
\end{subequations}
where $v_0=u_0-d_{1}(0)g(x)\in \mathbb{H}^2_{(0)}$ since $u_0\in \mathbb{H}^2_{(0)}$.
\vspace*{6pt}
\newline
Note that $|g(x)|\leq |b_1|+|b_2|+|c_1|+|c_2|:=g_0$, $|g_x|\leq |b_2|+2|c_1|+3|c_2|:=g_1$ and $|g_{xx}|\leq 2|c_1|+6|c_2|:=g_2$. Let $\tilde{\rho}(t,r)=g_0|d_{1t}(t)|+\mu g_2|d_1(t)|+\textcolor{black}{2^\gamma(1+g_1^\gamma|d_1(t)|^\gamma})\times\rho(t,r+2g_0|d_1(t)|)$, which is continuous in $t$ and $r$. \textcolor{black}{One may verify that $\tilde{f}(t,x,v,p)$ satisfies the structural conditions \eqref{AF} with $ \tilde{\rho}(t,r)$ instead of $\rho(t,r)$}.
According to Proposition \ref{Proposition 1}, \eqref{Equ. 4} has a unique solution $v\in  C^{0}(\mathbb{R}_{\geq 0}\times (0,1);\mathbb{H}^2_{(0)})\cap C^{1}(\mathbb{R}_{+}\times (0,1);L^2(0,1))$. Finally, $u=v+d_1(t)g(x)$ is the unique solution of \eqref{Equ. 1} and satisfies \eqref{B. 1}, \eqref{B. 2} and \eqref{B. 3}.
\vspace*{6pt}
\newline
For $b_1^2+b_2^2\neq 1$, we set
$\tilde{b}_i=\frac{b_i}{\sqrt{b_1^2+b_2^2}} (i=1,2),\ \tilde{d}_1=\frac{d_1}{\sqrt{b_1^2+b_2^2}}.$ Then the boundary condition \eqref{B. 2} is equivalent to $ \tilde{b}_1u(t,0)+\tilde{b}_2u_x(t,0)=\tilde{d}_{1}(t),$ where $\tilde{b}_1^2+\tilde{b}_2^2= 1 $. Therefore, \eqref{Equ. 1} has a unique solution $u\in  C^{0}(\mathbb{R}_{\geq 0}\times (0,1);\mathbb{H}^2_{(0)})\cap C^{1}(\mathbb{R}_{+}\times (0,1);L^2(0,1))$.
\hfill $\blacksquare$
\end{pf*}
\section{Stability Assessment}\label{Sec: ISS}
In stability analysis, we choose the energy of the system, $E(t) = \|u(t,\cdot)\|^2$, as the Lyapunov functional candidate. Let $\mathbb{H}^2_{(0)}$ be defined as in Section \ref{Sec: Well-posedness}. Note that in order to apply Lemma~\ref{Lemma 3} and Lemma~\ref{Lemma 2} to deal with the terms of $u_x$ on the boundaries, we always assume that $b_2\neq 0$ \textcolor{black}{(i.e., we consider the problem with Robin (or Neumann) boundary conditions)}.

\subsection{The case where the function $f(t,x,u,p)$ is in a general form}
We assume that \textcolor{black}{there exists \textcolor{black}{$d(t,x) \in C^1(\mathbb{R}_{\geq 0}\times (0,1); \mathbb{R})$ such that}}
\begin{align*}
 f(t,x,u,p)u\leq M_1u^2+(|d(t,x)|+M_2|p|)|u|,
\end{align*}
for a.e. $x\in (0,1)$ and all $ t\in\mathbb{R}_{\geq 0} , u\in\mathbb{R},p\in\mathbb{R}$, and {$M_1\in\mathbb{R}$ and $M_2\in\mathbb{R}_{\geq 0}$} are constants.  \textcolor{black}{Note that $d(t,x)$ can be used to describe the disturbance in the domain.} \textcolor{black}{For simplicity, we assume that $|d(x,t)|\leq |d_2(t)|$ for almost all $x\in (0,1)$ and any $t>0$, where $d_2\in C^1(\mathbb{R}_{\geq 0}; \mathbb{R})$, i.e., we assume that}
\begin{align}
 f(t,x,u,p)u\leq M_1u^2+(|d_2(t)|+M_2|p|)|u|,\label{AF-1-2}
\end{align}
for a.e. $x\in (0,1)$ and all $ t\in\mathbb{R}_{\geq 0} , u\in\mathbb{R},p\in\mathbb{R}$.
\textcolor{black}{\begin{definition}
System~\eqref{Equ. 1} with ~\eqref{B} is said to be input-to-state stable (ISS), or respectively integral input-to-state stable (iISS), w.r.t. the disturbances $d_1(t)$ and $d_2(t)$, if there exist functions $\beta\in \mathcal {K}\mathcal {L},\theta_1,\theta_2\in \mathcal {K}_{\infty}$ and $ \gamma_1, \gamma_2,\in \mathcal {K}$ such that the solution of \eqref{Equ. 1} with ~\eqref{B} satisfies
\begin{align}\label{Eq: ISS def}
\begin{split}
    \|u(t,\cdot)\|\leq &\beta( \|{u_0}\|,t)
      +\gamma_1(\|d_1\|_{L^{\infty}(0,t)}) \\
      &+\gamma_2(\|d_2\|_{L^{\infty}(0,t)}),\ \forall t\geq 0,
\end{split}
\end{align}
or respectively
\begin{align}\label{Eq: ISS def2}
\begin{split}
\|u(t,\cdot)\|\leq & \beta( \|{u_0}\|,t)+\theta_1\bigg(\int_{0}^{t} \gamma_1(|d_1(s)|)\bigg)\\
      &+\theta_2\bigg(\int_{0}^{t} \gamma_2(|d_2(s)|)\bigg),\ \forall t\geq 0.
\end{split}
\end{align}
Moreover, System~\eqref{Equ. 1} with ~\eqref{B} is said to be exponential input-to-state stable (EISS), or exponential integral input-to-state stable (EiISS), w.r.t. the disturbances $d_1(t)$ and $d_2(t)$, if there exist $\beta'\in \mathcal {K}_{\infty}$ and a constat $\lambda > 0$ such that  $\beta( \|{u_0}\|,t) \leq \beta'(\|{u_0}\|)e^{-\lambda t}$ in \eqref{Eq: ISS def} or \eqref{Eq: ISS def2}.
\end{definition}}
\textcolor{black}{
\begin{remark}\label{Rem: ISS}
While the ISS typically refers to norm-estimates for the input/disturbance in the $L^\infty$-norm, other norms can also be considered. It should be mentioned that the latter case usually relates to the integration of the input/disturbance and can be defined as the ``integral input-to-state stability (iISS)" (see, e.g., \cite[Definition~2.6]{Jacob:2016}). This property differs from the ISS in the sense that it allows for unbounded inputs that have ``finite energy" \cite{Sontag:1998}. There indeed exist many practically relevant systems that are iISS, but not ISS (see, e.g., \cite{karafyllis2011stability,Mironchenko:2014a} for more detailed discussions).
\end{remark}}
In order to obtain the stability of the system, we need some additional assumptions on $a_1,a_2,b_1,b_2,M_1$ and $M_2$. Specifically, if $a_2\neq 0$, we make the following assumptions.
\begin{assumption}\label{(A1.1)} Suppose that \eqref{+ab1} holds. Moreover, suppose that there \textcolor{black}{exist} $A'_1,A'_2,A'_3\in \mathbb{R}_{\geq 0}$ satisfying $A'_1+A'_2+A'_3=\mu $ and
\begin{align}\label{ab11}
 \frac{a_1}{a_2}\mu\geq 2A'_2,\ \frac{b_1}{b_2}\mu< A'_1.
\end{align}
Assume further that there exists $\varepsilon_0 \in \mathbb{R}_{+}$ such that
\begin{align}\label{ab12}
 M_1+\frac{\varepsilon_0M_2}{2}<A'_2-2A'_1,\  \frac{M_2}{2\varepsilon_0}\leq A'_3.
\end{align}
\end{assumption}
\begin{assumption}\label{(A1.2)} Suppose that \eqref{+ab2} holds. Moreover, suppose there \textcolor{black}{exist} $B'_1,B'_2,B'_3\in \mathbb{R}_{\geq 0}$ satisfying $B'_1+B'_2+B'_3=\mu $ and
\begin{align}\label{ab13}
 \frac{a_1}{a_2}\mu\geq -B'_2,\ \frac{b_1}{b_2}\mu< -2B'_1.
\end{align}
Assume further that there exists $\varepsilon_0 \in \mathbb{R}_{+}$ such that
\begin{align}\label{ab14}
 M_1+\frac{\varepsilon_0M_2}{2}<B'_1-2B'_2,\ \frac{M_2}{2\varepsilon_0}\leq B'_3.
\end{align}
\end{assumption}
If $a_2=0$, we make the following assumption.
\begin{assumption}\label{(A2.1)} We assume that $\frac{b_1}{b_2}< \frac{1}{2} $, which guarantees that there \textcolor{black}{exist} $A'_1,A'_2,A'_3\in \mathbb{R}_{\geq 0}$ satisfying $A'_1+A'_2+A'_3=\mu $ and $
 \frac{b_1}{b_2}\mu< A'_1.
$
Assume further that there exists $\varepsilon_0 \in \mathbb{R}_{+}$ such that
\begin{align*}
 M_1+\frac{\varepsilon_0M_2}{2}<A'_2-2A'_1,\ \frac{M_2}{2\varepsilon_0}\leq A'_3.
\end{align*}
\end{assumption}
%
\begin{theorem}\label{Theorem 8}
Let $u\in  C^{0}(\mathbb{R}_{\geq 0}\times (0,1);\mathbb{H}^2_{(0)})\cap C^{1}(\mathbb{R}_{+}\times (0,1);L^2(0,1))$ be the unique solution of \eqref{Equ. 1}, \eqref{B. 1},\eqref{B. 2} and \eqref{B. 3}. Under Assumption~\ref{(A1.1)}, or Assumption~\ref{(A1.2)}, or Assumption~\ref{(A2.1)}, 
\textcolor{black}{System~\eqref{Equ. 1} with ~\eqref{B} is EiISS and EISS having the estimates:}
\begin{align}\label{Eq: ISS 1}
\|u(t,\cdot)\|^2\leq &\|u(0,\cdot)\|^2e^{-C_0t}+C_1\int_{0}^t |d_1(s)|^2\text{d}s \notag\\
 & \;\;\;+ C_2\int_{0}^t|d_2(s)|^2\text{d}s,
\end{align}
and
\begin{align}\label{Eq: ISS 2}
\|u(t,\cdot)\|^2\leq &\|u(0,\cdot)\|^2e^{-C_0t}+\Big(1-e^{-C_0t}\Big) \notag\\
   & \times\Big(C_1\|d_1\|_{L^{\infty}(0,t)}^2+C_2\|d_2\|_{L^{\infty}(0,t)}^2\Big)
\end{align}
for some positive constants $C_0,C_1,C_2$.
\end{theorem}
\begin{pf*}{Proof.} We prove first the case for $a_2\neq 0$ under Assumption~\ref{(A1.1)}.
Multiplying \eqref{Equ. 1} with $u$ and integrating over $[0,1]$, we have
\begin{align*}
\int_{0}^1u_tu\text{d}x&-\mu \int_{0}^1u_{xx}u\text{d}x
=  \int_{0}^1f(t,x,u,u_x)u\text{d}x\notag\\
\leq & \int_{0}^1\big((|d_2(t)|+M_2|u_x|)|u|+M_1u^2\big)\text{d}x := I_1,
\end{align*}
which is
\begin{align*}
\frac{\text{d}}{\text{d}t}\|u\|^2-\mu u_{x}(t,1)u(t,1)+\mu u_{x}(t,0)u(t,0)+\mu \|u_x\|^2\leq I_1.
\end{align*}
By \eqref{B. 1}, \eqref{B. 2} and Young's inequality, it follows
\begin{align}\label{62810}
&\frac{\text{d}}{\text{d}t}\|u\|^2+\mu \|u_x\|^2\notag\\
\leq & I_1-\frac{1}{b_2}d_1(t)\mu u(t,0)+\frac{b_1}{b_2}\mu u^2(t,0)-\frac{a_1}{a_2}\mu u^2(t,1)\notag\\
\leq & I_1+\frac{\mu |d_1(t)|^2}{2\varepsilon_1b_2^2}
+\bigg(\frac{b_1}{b_2}+\frac{\varepsilon_1}{2}\bigg)\mu u^2(t,0)-\frac{a_1}{a_2}\mu u^2(t,1).
\end{align}
By Young's inequality, we have
\begin{align}\label{62811}
I_1\leq
%
\frac{|d_2(t)|^2}{2\varepsilon_2}+\bigg(\frac{\varepsilon_2}{2}+M_1
+\frac{\varepsilon_3M_2}{2}\bigg)\|u\|^2+\frac{M_2}{2\varepsilon_3}\|u_x\|^2.
\end{align}
Then we infer from \eqref{62810}, \eqref{62811} and \eqref{62812} that
\begin{align}\label{62813}
\frac{\text{d}}{\text{d}t}\|u\|^2&+A'_1u^2(t,0)+(A'_2-2A'_1)\|u\|^2-2A'_2u^2(t,1)\notag\\
 &+A'_3\|u_x\|^2\notag\\
=&\frac{\text{d}}{\text{d}t}\|u\|^2+A'_1(u^2(t,0)-2\|u\|^2)\notag\\
 &+A'_2(\|u\|^2-2u^2(t,1) )+A'_3\|u_x\|^2\notag\\
\leq&\frac{\text{d}}{\text{d}t}\|u\|^2+A'_1\|u_x\|^2+A'_2\|u_x\|^2+A'_3\|u_x\|^2\notag\\
=&
\frac{\text{d}}{\text{d}t}\|u\|^2+\mu \|u_x\|^2\notag\\
\leq &\frac{\mu |d_1(t)|^2}{2\varepsilon_1b_2^2}+\frac{|d_2(t)|^2}{2\varepsilon_2}+\bigg(\frac{b_1}{b_2}
+\frac{\varepsilon_1}{2}\bigg)\mu u^2(t,0)\notag\\
&-\frac{a_1}{a_2}\mu u^2(t,1) +\bigg(\frac{\varepsilon_2}{2}+M_1
 +\frac{\varepsilon_3M_2}{2}\bigg)\|u\|^2 \notag\\
&+\frac{M_2}{2\varepsilon_3}\|u_x\|^2.
\end{align}
Recalling \eqref{ab11} and \eqref{ab12}, one may choose $\varepsilon_3=\varepsilon_0$ and $\varepsilon_1, \varepsilon_2>0$ small enough such that
\begin{align*}
&C_0:= A'_2-2A'_1-\bigg(\frac{\varepsilon_2}{2}+M_1
+\frac{\varepsilon_3M_2}{2}\bigg)>0,\\
& \bigg(\frac{b_1}{b_2}+\frac{\varepsilon_1}{2}\bigg)\mu \leq A'_1,\ \frac{M_2}{2\varepsilon_3}\leq A'_3.
\end{align*}
Then we have
\begin{align}
\frac{\text{d}}{\text{d}t}\|u\|^2&\leq - C_0\|u\|^2+\frac{\textcolor{black}{\mu|d_1(t)|^2}}{2\varepsilon_1b_2^2}+\frac{|d_2(t)|^2}{2\varepsilon_2}\notag\\
&:=-C_0\|u\|^2+C_1|d_1(t)|^2+C_2|d_2(t)|^2\label{+18}\\
&\leq -C_0\|u\|^2+C_1\|d_1\|_{L^{\infty}(0,t)}^2
+C_2\|d_2\|_{L^{\infty}(0,t)}^2.\label{+19}
\end{align}
By \eqref{+18} and Gronwall's inequality, we obtain \eqref{Eq: ISS 1}.
By \eqref{+19} and Gronwall's inequality, we obtain \eqref{Eq: ISS 2}.
\vspace*{6pt}
\newline
For $a_2\neq 0$ and under Assumption~\ref{(A1.2)}, it suffices to note that by Lemma~\ref{Lemma 3} and Lemma~\ref{Lemma 2}, we have $\|u_x\|^2\geq u^2(t,1)-2\|u\|^2$  and  $\|u_x\|^2\geq \|u\|^2-2u^2(t,0)$. Then proceeding as above, one may get
\begin{align*}
\frac{\text{d}}{\text{d}t}\|u\|^2&+(B'_1-2B'_2)\|u\|^2+B'_3\|u_x\|^2  -2B'_1 u^2(t,0) \notag\\
                   & +B'_2u^2(t,1) \notag\\
\leq & \frac{\mu|d_1(t)|^2}{2\varepsilon_1b_2^2}+\frac{|d_2(t)|^2}{2\varepsilon_2}+\bigg(\frac{b_1}{b_2}
+\frac{\varepsilon_1}{2}\bigg)\mu u^2(t,0)\notag\\
 &-\frac{a_1}{a_2}\mu u^2(t,1) +\bigg(\frac{\varepsilon_2}{2}+M_1
+\frac{\varepsilon_3M_2}{2}\bigg)\|u\|^2 \notag\\
&+\frac{M_2}{2\varepsilon_3}\|u_x\|^2.
\end{align*}
\textcolor{black}{The ISS can be established as well.}
\vspace*{6pt}
\newline
Now for $a_2=0$, it suffices to note that $u(t,1)=0$. Under Assumption~\ref{(A2.1)}, the ISS can be obtained as above.
\hfill $\blacksquare$
\end{pf*}
\begin{remark} \label{Remark 2}
It should be noticed that the assumptions that \eqref{+ab1} (or \eqref{+ab2}) holds in Assumption~\ref{(A1.1)} (or Assumption~\ref{(A1.2)}) are only for assuring the existence of a solution. For ISS assessment, it suffices to relax these assumptions to \eqref{ab11} and \eqref{ab12} (or \eqref{ab13} and \eqref{ab14}), or some other weaker conditions.
\end{remark}
\subsection{The case where the function $f(t,x,u,p)$ has a special form}
In the following part, we assume that $f(t,x,u,p)$ is with the form
\begin{align}
 f(t,x,u,p)=\textcolor{black}{d(t,x)}+ M_1u+M_2p,\label{AF-1-3}
\end{align}
where $M_1, M_2\in \mathbb{R}$ are constants, $|d(t,x)|\leq |d_2(t)|$ for a.e. $x\in (0,1)$, $d_2\in C(\mathbb{R}_{\geq 0}; \mathbb{R})$.
\textcolor{black}{As $f(t,x,u,p)$ grows lineally w.r.t. $u$ and $p$, the conditions given in Assumption~\ref{(A1.1)}, Assumption~\ref{(A1.2)}, and Assumption~\ref{(A2.1)} can be relaxed.}

For the case where $a_2\neq 0$, we make the following assumptions.
\begin{assumption}\label{(A3a)}  Suppose that there \textcolor{black}{exist} $A'_1,A'_2\in \mathbb{R}_{\geq 0}$ satisfying $A'_1+A'_2=\mu $ and
\begin{align*}
- \frac{a_1}{a_2}\mu+\frac{M_2}{2}\leq -2A'_2,\ \frac{b_1}{b_2}\mu-\frac{M_2}{2}< A'_1,\ M_1<A'_2-2A'_1.
\end{align*}
\end{assumption}
\begin{assumption}\label{(A3b)} Suppose that there \textcolor{black}{exist} $B'_1,B'_2\in \mathbb{R}_{\geq 0}$ satisfying $B'_1+B'_2=\mu $ and
\begin{align*}
 -\frac{a_1}{a_2}\mu+\frac{M_2}{2}\leq B'_2,\ \frac{b_1}{b_2}\mu-\frac{M_2}{2}< -2B'_1,\ M_1<B'_1-2B'_2.
\end{align*}
\end{assumption}

For the case where $a_2= 0$, we make the following assumption.
\begin{assumption}\label{(A3c)} We assume that there \textcolor{black}{exist} $A'_1,A'_2\in \mathbb{R}_{\geq 0}$ satisfying $A'_1+A'_2=\mu $ such that
\begin{align*}
 \frac{b_1}{b_2}\mu-\frac{M_2}{2}<A_1',\ M_1<A'_2-2A'_1.
\end{align*}
\end{assumption}

\begin{theorem}\label{Theorem 9}
Let $u\in  C^{0}(\mathbb{R}_{\geq 0}\times (0,1);\mathbb{H}^2_{(0)})\cap C^{1}(\mathbb{R}_{+}\times (0,1);L^2(0,1))$ be the unique solution of \eqref{Equ. 1}, \eqref{B. 1},\eqref{B. 2} and \eqref{B. 3}. Under Assumption~\ref{(A3a)}, or Assumption~\ref{(A3b)}, or Assumption~\ref{(A3c)}, \textcolor{black}{System~\eqref{Equ. 1} with ~\eqref{B} is EiISS and EISS having the estimates:}
\begin{align*}
\|u(t,\cdot)\|^2\leq &\|u(0,\cdot)\|^2e^{-C_3t}+C_4\int_{0}^t |d_1(s)|^2\text{d}s \notag\\
 & \;\;\;+ C_5\int_{0}^t|d_2(s)|^2\text{d}s,
\end{align*}
and
\begin{align*}
\|u(t,\cdot)\|^2\leq &\|u(0,\cdot)\|^2e^{-C_3t}+ \Big(1-e^{-C_0t}\Big) \notag\\
   & \;\;\; \times\Big(C_4\|d_1\|_{L^{\infty}(0,t)}^2+C_5\|d_2\|_{L^{\infty}(0,t)}^2\Big)
\end{align*}
for some positive constants $C_3,C_4,C_5$.
\end{theorem}
\begin{pf*}{Proof.} We proceed as in Theorem \ref{Theorem 8} and only prove the result under  Assumption \ref{(A3a)}. Multiplying \eqref{Equ. 1} with $u$ and integrating over $[0,1]$, we have
\begin{align*}
\frac{\text{d}}{\text{d}t}\|u\|^2&-\mu u_{x}(t,1)u(t,1)+\mu u_{x}(t,0)u(t,0)+\mu \|u_x\|^2\notag\\
= & \int_{0}^1f(t,x,u,u_x)u\text{d}x\notag\\
\leq & \int_{0}^1\big(|d_2(t)||u|+M_1u^2+M_2u_xu\big)\text{d}x := I_2.
\end{align*}
Note that $\int_{0}^1u_xu\text{d}x=\frac{1}{2}u^2(t,x)|^{x=1}_{x=0}=\frac{1}{2}(u^2(t,1)- u^2(t,0))$. It follows
\begin{align*}
I_2\leq& \frac{|d_2(t)|^2}{2\varepsilon_2}+\frac{\varepsilon_2}{2}\|u\|^2
+M_1\|u\|^2+\frac{M_2}{2}(u^2(t,1)- u^2(t,0))\notag\\
=&\frac{|d_2(t)|^2}{2\varepsilon_2}+\bigg(\frac{\varepsilon_2}{2}+M_1
\bigg)\|u\|^2+\frac{M_2}{2}(u^2(t,1)- u^2(t,0)).
\end{align*}
Then we have
\begin{align*}
&\frac{\text{d}}{\text{d}t}\|u\|^2+\mu \|u_x\|^2\notag\\
\leq & I_2-\frac{1}{b_2}d_1(t)\mu u(t,0)+\frac{b_1}{b_2}\mu u^2(t,0)-\frac{a_1}{a_2}\mu u^2(t,1)\notag\\
\leq & \frac{\mu |d_1(t)|^2}{2\varepsilon_1b_2^2}+\frac{|d_2(t)|^2}{2\varepsilon_2}+\bigg(\frac{\varepsilon_2}{2}+M_1
\bigg)\|u\|^2\notag\\
&+\bigg(\frac{b_1\mu}{b_2}-\frac{M_2}{2}+\frac{\varepsilon_1\mu}{2}\bigg) u^2(t,0)+\bigg(\frac{M_2}{2}-\frac{a_1}{a_2}\mu \bigg)u^2(t,1).
\end{align*}
We get by splitting $\mu \|u_x\|^2$ as in \eqref{62813} and using \eqref{62812}
\begin{align*}
\frac{\text{d}}{\text{d}t}\|u\|^2&+A'_1u^2(t,0)+(A'_2-2A'_1)\|u\|^2-2A'_2u^2(t,1)
 \notag\\
=&\frac{\text{d}}{\text{d}t}\|u\|^2+A'_1(u^2(t,0)-2\|u\|^2)\notag\\
 &+A'_2(\|u\|^2-2u^2(t,1) )\notag\\
\leq&
\frac{\text{d}}{\text{d}t}\|u\|^2+\mu \|u_x\|^2\notag\\
\leq &\frac{\mu |d_1(t)|^2}{2\varepsilon_1b_2^2}+\frac{|d_2(t)|^2}{2\varepsilon_2}+\bigg(\frac{\varepsilon_2}{2}+M_1
\bigg)\|u\|^2\notag\\
&+\bigg(\frac{b_1\mu}{b_2}-\frac{M_2}{2}+\frac{\varepsilon_1\mu}{2}\bigg) u^2(t,0)\notag\\
&+\bigg(\frac{M_2}{2}-\frac{a_1}{a_2}\mu \bigg)u^2(t,1).
\end{align*}
Choosing $ \varepsilon_1,\varepsilon_2$ small enough, such that
\begin{align*}
&\frac{b_1\mu}{b_2}-\frac{M_2}{2}+\frac{\varepsilon_1\mu}{2}\leq A_1',\notag\\
&\frac{\varepsilon_2}{2}+M_1
<A_2'-2A_1',\frac{M_2}{2}-\frac{a_1}{a_2}\mu \leq -2A_2'.
\end{align*}
Then we have
\begin{align*}
\frac{\text{d}}{\text{d}t}\|u\|^2&\leq -C_3\|u\|^2+C_4|d_1(t)|^2+C_5|d_2(t)|^2\\
&\leq -C_3\|u\|^2+C_4\|d_1\|_{L^{\infty}(0,t)}^2
+C_5\|d_2\|_{L^{\infty}(0,t)}^2.
\end{align*}
Finally, one may obtain the desired results by Gronwall's inequality.
\hfill $\blacksquare$
\end{pf*}
\begin{remark}
\textcolor{black}{Note that with Assumption~\ref{(A3a)}, or Assumption~\ref{(A3b)}, or Assumption~\ref{(A3c)}, Proposition~\ref{Proposition 1} guarantees that the operator $\mathcal{A}$ generates an exponentially stable semi-group, when $d_1(t)=0$. Then, it follows directly from Proposition~4 in \cite{Mironchenko:2014a} or Proposition~2.13 in \cite{Jacob:2016} that System~\eqref{Equ. 1} with ~\eqref{B} is ISS w.r.t. $d_2$ in $L^p$-norm ($p \geq  1$). For $d_1(t)\neq 0$, the ISS w.r.t. disturbances in $L^\infty$-norm is obtained in \cite{Jacob:2016,karafyllis2017siam}. This result is weaker than that obtained in Theorem~\ref{Theorem 9}, which gives an ISS w.r.t. disturbances in $L^2$-norm. It should be mentioned that a strict Lyapunov functional has also been constructed to establish the ISS w.r.t. in-domain disturbances for a semi-linear parabolic PDE with periodic boundary conditions \cite[Theorem~3]{Mazenc:2011}.
}
\end{remark}

\section{Illustration Examples}\label{Sec: Examples}
Two examples are used to illustrate the developed results.
\subsection{Ginzburg-Landau equations with real coefficients}
Consider first the Ginzburg-Landau equation with real coefficients (see, e.g., \cite{Kogelman:1970})
\begin{align}
u_t= {\mu u_{xx}+\alpha u}-\beta |u|^2u,\label{Landau 1}
\end{align}
and the generalized Ginzburg-Landau equation with real coefficients (see, e.g., \cite{Guo:1994})
\begin{align}
u_t= {\mu u_{xx}+\alpha u}-\beta |u|^2u-\gamma |u|^4u+\lambda u_x,\label{Landau 2}
\end{align}
under the Robin (or Neumann) boundary conditions
\begin{align*}
      u(x,1)=0,          b_1u(x,0)+b_2u_x(t,0)&=d(t),
      \end{align*}
where $\mu ,\beta,\gamma>0,\alpha,\lambda,b_1,b_2\in\mathbb{R},b_2\neq 0$ and $d\in C^{2}(\mathbb{R}_{\geq 0};\mathbb{R})$.

In the above boundary conditions, $a_1=1$, $a_2=0$, $d_1(t)=d(t)$, and $d_2(t)=0$.
In \eqref{Landau 1}, $f(t,x,u,p)=\alpha u-\beta |u|^2u$. In \eqref{Landau 2}, $f(t,x,u,p)=\alpha u-\beta |u|^2u-\gamma |u|^4u+\lambda p$. \textcolor{black}{In both cases, $f(t,x,u,p)$ satisfies the structural conditions} \eqref{AF}. Assume that $\frac{b_1}{b_2}\leq \frac{1}{2}$, then there exists a unique real solution of \eqref{Landau 1} and \eqref{Landau 2} respectively. \newline
Now for \eqref{Landau 1}, $f(t,x,u,p)$ satisfies \eqref{AF-1-2} with
$
M_1=\alpha,M_2=0.
$
If we assume further that
$
  \frac{b_1}{b_2}< \frac{1}{3},\alpha<0,
$
and set $A_1'=\frac{1}{3}\mu,A_2'=\frac{2}{3}\mu, {A_3'=0}$, then Assumption \ref{(A2.1)} holds. Therefore \eqref{Landau 1} is ISS.
\newline
For \eqref{Landau 2}, $f(t,x,u,p)$ satisfies \eqref{AF-1-2} with
$
M_1=\alpha,M_2=|\lambda|.
$
If we assume further that
$
  \frac{b_1}{b_2}< \frac{1}{4},\alpha+{|\lambda|}<0,|\lambda|\leq \mu,
$
and set $A_1'=\frac{1}{4}\mu,A_2'=\frac{1}{2}\mu, {A_3'=\frac{1}{4}\mu},\varepsilon_0=2$, then Assumption \ref{(A2.1)} holds. Therefore \eqref{Landau 2} is ISS.

\subsection{1-$D$ transport partial differential equation}
We consider the following 1-$D$ transport PDE:
\begin{align}\label{Example 1}
u_t=\mu u_{xx}-mu_{x}-nu,
\end{align}
under the following boundary conditions
\begin{align}
\begin{split}\label{Eq: example 2 BC}
  u_x(t,1)&=\left(\frac{m}{2\mu}-a\right)u(t,1),\\
  u_x(t,0)&=\left(\frac{m}{2\mu}-b\right)u(t,0)+d(t),
\end{split}
\end{align}
where $\mu>0,m\geq 0, {n,a,b}\in\mathbb{R}$ and $d\in C^{2}(\mathbb{R}_{\geq 0};\mathbb{R})$.
\vspace*{6pt}
\newline
In order to make the manipulations easier,
we set $w(t,x)=e^{\frac{mx}{2\mu}}u(t,x)$. We can then transform the PDE \eqref{Example 1} with boundary conditions to the following problem (see also \cite{Karafyllis:2016}):
\begin{subequations}\label{Example 2.2}
\begin{align}
&w_t=\mu w_{xx}-\bigg(\frac{m^2}{4\mu}+n\bigg)w,\\
&w_x(t,1)=-aw(t,1),\\
&w_x(t,0)=-bw(t,0)+d(t).
\end{align}
\end{subequations}
In this case,
\begin{align*}
&a_1=a,a_2=1,b_1=b,b_2=1,\notag\\
&d_1(t)=d(t),d_2(t)=0,M_1=-\bigg(\frac{m^2}{4\mu}+n\bigg),M_2=0.
\end{align*}
If we assume that
$
a\geq\frac{4}{3},b<\frac{1}{3},\frac{m^2}{4\mu}+n>0,
$
and set $ A_1=\frac{1}{3},A_2=\frac{2}{3},A_1'=\frac{1}{3}\mu,A_2'=\frac{2}{3}\mu$, then condition \eqref{+ab1} and Assumption \ref{(A3a)} hold.
If we assume that
$
a\geq-\frac{1}{3},b<-\frac{4}{3},\frac{m^2}{4\mu}+n>0,
$
and set $ B_1=\frac{2}{3},B_2=\frac{1}{3},B_1'=\frac{2}{3}\mu,B_2'=\frac{1}{3}\mu$, then condition \eqref{+ab2} and Assumption \ref{(A3b)} hold. Under the above two assumptions, \eqref{Example 2.2} has a unique solution and \eqref{Example 2.2} is ISS, and so is \eqref{Example 1}.
%
%
\begin{remark} If it is easy to fix $A'_i, B'_i\ (i=1,2)$, one may verify the conditions in Assumption~\ref{(A3a)} and Assumption~\ref{(A3b)} to conclude the ISS of \eqref{Example 1} directly. It should be noticed that Assumption~\ref{(A3a)} or Assumption~\ref{(A3b)} is not a necessary condition for the ISS. Therefore, the system may be ISS even if Assumption~\ref{(A3a)} or Assumption~\ref{(A3b)} fails (see also Remark \ref{Remark 2}).
\end{remark}
\begin{remark}
The system \eqref{Example 1} was considered in \cite{Karafyllis:2016} under the boundary conditions:
\begin{enumerate}[(i)]
\item Dirichlet boundary conditions:\\
      $
       u(t,1)=0, \; u(t,0)=d(t).
      $
\item Robin (or Neumann) boundary conditions:\\
      $
      u(t,0)=d(t),u_x(t,1)=\left(\frac{m}{2\mu}-a\right)u(t,1), (a\geq 0).
     $
\end{enumerate}
In the above boundary conditions, $b_2=0$, which is slightly different from \eqref{Eq: example 2 BC}. The ISS property of \eqref{Example 1} was obtained by Parseval's identity and the expansions of eigenfunctions of the Sturm-Liouville operator under the same assumption that $\frac{m^2}{4\mu}+n>0$. \textcolor{black}{For the system \eqref{Example 1} ($m=0$) with Dirichlet boundary conditions under a boundary state feedback, the ISS in $L^p$-norm ($p\in (2,+\infty)$) is established in \cite{Mironchenko:2017} by the monotonicity-based method.}
\end{remark}

\section{Conclusion}\label{Sec: Conclusion}
This paper demonstrated via the considered semi-linear PDE that the ISS property with respect to \textcolor{black}{Robin (or Neumann) boundary} disturbances can be derived from suitable Lyapunov functionals. The obtained results confirmed that the appearance of the derivatives of boundary disturbances in the ISS estimates can be avoided by directly dealing with the boundary conditions with disturbances. Compared to the work reported in \cite{Karafyllis:2016,Karafyllis:2016a,karafyllis2017siam}, the application of Lyapunov functionals in the establishment of \emph{a priori} estimates of the solution seems to be less computationally demanding. Therefore, it can be expected that the developed techniques may be applicable in the study of ISS properties for a wider class of PDEs. \textcolor{black}{Finally, it should be mentioned that the technique developed in this work cannot deal with the ISS w.r.t. Dirichlet boundary disturbances, which is the case where $b_2 = 0$ in \eqref{B}. To tackle this type of problems, a method is developed in a parallel work \cite{Zheng:2017}.}

\bibliographystyle{plain}        
\bibliography{References}        


\end{document}